\documentclass[12pt,a4paper]{amsart}
\usepackage{a4}
\usepackage[english]{babel} 
\usepackage{amsmath,amssymb,amsthm}
\def\leq {\leqslant}
\def\le {\leqslant}
\def\ge {\geqslant}
\def\geq {\geqslant}
\def\@bibitem[#1]#2{\item\@biblabel{#1}.\if@filesw
{\def\protect##1{\string##1\space}\immediate\write
\@auxout{\string\bibcite{#2}{#1}}}\fi\ignorespaces\@showtag{#2}}

\theoremstyle{plain}
\newtheorem{theorem}{Theorem}[section]
\newtheorem{rem}{Remark}[section]
\newtheorem{lemma}{Lemma}[section]
\newtheorem{op}{Definition}[section]

\renewcommand{\theequation}%
{\arabic{section}.\arabic{equation}}
\begin{document}

\title{On the exactness of the conditions of embedding theorems for spaces of functions with mixed logarithmic smoothness}
\author{ G. Akishev}
\address{ Lomonosov Moscow University, Kazakhstan Branch \\
Str. Kazhymukan, 11 \\
010010, Astana, Kazakhstan}

\address{
Institute of mathematics and mathematical modeling,\\
Pushkin str, 125, \\
050010, Almaty, \\
 Kazakhstan
 }

\maketitle

\begin{quote}
\noindent{\bf Abstract. }
The article considers the Lorentz space $L_{p,\tau}(\mathbb{T}^{m})$, $2\pi$ of periodic functions of many variables and $S_{p,\tau,\theta}^{0, \overline{b}}\mathbf{B}$, $S_{p, \tau, \theta}^{0, \overline{b}}B$ --- spaces of functions with mixed logarithmic smoothness. The article establishes necessary and sufficient conditions for embedding the spaces $S_{p, \tau, \theta}^{0, \overline{b}}\mathbf{B}$ and $S_{p, \tau, \theta}^{ 0, \overline{b}}B$ into each other.
\end{quote}

\vspace*{0.2 cm}


\vspace*{0.2 cm}

\section*{Introduction}
\label{intro} 

Let $\mathbb{N}$, $\mathbb{Z}$, $\mathbb{R}$ be sets of natural, integer, and real numbers, respectively, and $\mathbb{Z}_{+} := \mathbb{N }\cup\{0\}$, and $\mathbb{R}^{m}$ be an $m$--dimensional Euclidean space of points $\overline{x}=(x_{1},\ldots,x_{m})$ with real coordinates;
$\mathbb{T}^{m}=[0, 2\pi)^{m}$ and $\mathbb{I}^{m}=[0, 1)^{m}$ be $m $--dimensional cubes. Further, $\mathbb{Z}^{m}$ and $\mathbb{Z}_{+}^{m}$ are $m$--fold Cartesian product of the sets $\mathrm{Z}$ and $\mathbb{Z }_{+}$, respectively.
 \smallskip

We denote by $L_{p,\tau}(\mathbb{T}^{m})$
the Lorentz space of all real-valued Lebesgue measurable functions $f$ that have $2\pi$--period in each variable and for which the quantity
\begin{equation*}
\|f\|_{p,\tau} = \left\{\frac{\tau}{p}\int\limits_{0}^{1}\biggl(f^{*}(t)
\biggr)^{\tau}t^{\frac{\tau}{p}-1}dt
\right\}^{\frac{1}{\tau}} , \,\, 1<
p<\infty, 1\leqslant \tau <\infty,
\end{equation*}
is finite, where $f^{*}(t)$ is a non-increasing rearrangement of the function $|f(2\pi\overline{x})|$, $\overline{x} \in \mathrm{I}^{m }$ (see \cite{1}, ch. 1, sec. 3).

In case when $\tau=p$, the Lorentz space $L_{p,\tau}(\mathbb{T}^{m})$ coincides with the Lebesgue space $L_{p}(\mathbb{T}^{m})$ with the norm $\|f\|_{p}=\|f\|_{p,p}$ (see, for example, \cite[ch. 1, sec. 1.1]{2}).

We denote by ${\mathring L}_{p, \tau}(\mathbb{T}^{m})$  the set of all functions $f\in
L_{p, \tau}(\mathbb{T}^{m})$,
such that
\begin{equation*}
\int\limits_{0}^{2\pi }f(\overline{x}) dx_{j}  =0,\;\;
j=1,...,m .
\end{equation*}

Below $a_{\overline{n}}(f)$ are Fourier coefficients of a function $f\in L_{1}\left(\mathrm{T}^{m} \right)$ by the system $\{e^{i\langle\overline{n}, 2\pi\overline{x}\rangle}\}_{\overline{n}\in\mathrm{Z}^{m}}$ and $\langle\overline{y},\overline{x}\rangle=\sum\limits_{j=1}^{m}y_{j} x_{j}$;
$$
	\delta_{\overline{s}}(f, 2\pi\overline{x})
	= \sum\limits_{\overline{n} \in \rho \left( \overline{s} \right)}a_{\overline{n} } (f) e^{i\langle\overline{n}, 2\pi\overline{x}\rangle },
$$
where
$$
	\rho (\overline{s})=\left\{\overline{k} =\left( k_{1}, \dots, k_{m} \right) \in \mathrm{Z}^{m}\colon [2^{s_{j} -1}] \leq \left| k_{j} \right|<2^{s_{j} } ,j=1, \dots, m\right\},
$$
and $[a]$ is an integer part of $a$, $\overline{s} = (s_{1}, \dots, s_{m}), s_{j} = 0, 1, 2,\ldots$.

The quantity
\begin{equation*}
Y_{l_{1},\ldots,l_{m}}(f)_{p, \tau} = \inf_{T_{l_{j}}} \|f-\sum_{j=1}^{m}T_{l_{j}}\|_{p, \tau}^{*} \;\;,\;\; l_{j} = 0,1,2,...
\end{equation*}
is called the best approximation of the ``angle`` of the function $f\in L_{p, \tau}(\mathbb{T}^{m})$ by trigonometric polynomials, where $T_{l_{j}} \in L_{p}(\mathbb{T}^{m})$ a trigonometric polynomial of order $l_{j}$ with respect to the variable $x_{j}, \,\, j = 1,\ldots,m$ (in the case of $\tau=p$, see \cite{3}--\cite{6}).

Throughout the paper, $A_{n}\asymp B_{n}$ means that there are positive numbers $C_{1}, C_{2}$ independent of $n\in \mathrm{N}$ such that $C_{1}A_{n}\leq B_{n} \leq C_{2}A_{n}$ for $n\in \mathrm{N}$ and $\log M$, where $\log M$ is the logarithm with base $2$ of the number $M>1$.

\begin{op}
The total smoothness modulus of the order $k$ of the function $f\in L_{p, \tau}(\mathbb{T}^{m})$ is the quantity (in the case of $\tau=p$, see  \cite{2}, ch. 4, sec. 4. 2)
\begin{equation*}
\omega_{k}(f, t)_{p, \tau} = \sup_{\|\overline{h}\|\leq t}\|\Delta_{\overline{h}}^{k}f\|_{p, \tau},
\end{equation*}
where $\Delta_{\overline{h}}^{k}f(\overline{x})$ is the total difference of order $k$, of the function $f$.
\end{op}

\begin{op}
The mixed smoothness modulus of a function of order $\overline{k}$ of the function $f\in L_{p, \tau}(\mathbb{T}^{m})$ is determined by the formula
\begin{equation*}
\omega_{\overline k}(f, \overline{t})_{p, \tau}=\omega_{k_{1},...k_{m}}(f, t_{1},..., t_{m})_{p, \tau}=\sup_{|h_{1}|\leq t_{1},...,|h_{m}|\leq t_{m}}\|\Delta_{\overline h}^{\overline k} (f)\|_{p, \tau},
\end{equation*}
where $\Delta_{\bar t}^{\bar k} f(\bar x) = \Delta_{t_{m}}^{k_{m}}(...
\Delta_{t_{1}}^{k_{1}}f(\overline{x}))$ --- mixed difference of order $\overline{k}$ in increments $\overline{t} = (t_{1},...,t_{m})$.
\end{op}

{\bf The space Besov} (\cite{7},  \cite{2}, ch. 4, sec. 4. 3).
Let $k \in \mathbb{N}$ and $k > r > 0$, $1\leq \theta \leq \infty$.
The Besov space $\mathbf{B}_{p, \theta}^{r}(\mathbb{T}^{m})$ is called the set of all functions $f\in L_{p}(\mathbb{T}^{m})$, $1\leq p < \infty$, for which
\begin{equation*}
\|f\|_{\mathbf{B}_{p, \theta}^{r}}: = \|f\|_{p} + \Bigl(\int_{0}^{1}(t^{-r}\omega_{k}(f, t)_{p})^{\theta} \frac{dt}{t}\Bigr)^{1/\theta} < \infty, \, \, 
\end{equation*}
for $1\leq \theta < \infty$ and 
\begin{equation*}
\|f\|_{\mathbf{B}_{p, \infty}^{r}}: =\|f\|_{p} +\sup_{t> o}t^{-r}\omega_{k}(f, t)_{p}< \infty
\end{equation*}
for $\theta = \infty$. 

The space $\mathbf{B}_{p, \infty}^{r}(\mathbb{T}^{m})$  
is denoted by the symbol $\mathbf{H}_{p}^{r}(\mathbb{T}^{m})$ (see \cite{2}).

The following generalization of the Besov space is defined in articles \cite{8}--\cite{10}, $\mathbf{B}_{p, \theta}^{r, b}(\mathbb{T}^{m})$ -- the set of all functions $f\in L_{p}(\mathbb{T}^{m})$, $1\leq p < \infty$, for which
\begin{equation*}
\|f\|_{\mathbf{B}_{p, \theta}^{r, b}}: = \|f\|_{p} + \Bigl(\int_{0}^{1}(t^{-r}(1 - \log t)^{b}\omega_{k}(f, t)_{p})^{\theta} \frac{dt}{t}\Bigr)^{1/\theta} < \infty, \, \, 
\end{equation*}
for $0< \theta \leq \infty$, $b> -1/\theta$, $k\in \mathbb{N}$, $k > r\ge 0$.

In the articles \cite{10}, \cite{11}, \cite{12} the space $B_{p, \theta}^{0, b}(\mathbb{T}^{m})$ is defined as the set of functions $f\in {\mathring L}_{p}(\mathbb{T}^{m})$ for which
\begin{equation*}
\left\{\sum\limits_{s \in \Bbb{Z}_{+}} (s + 1)^{b\theta} \|\sigma_{s}(f)\|_{p, \tau}^{\theta}\right\}^{\frac{1}{\theta}} < \infty,
\end{equation*}
for $1< p < \infty$, $b > -1/\theta$, $0< \theta < \infty$.

Relations between the spaces $B_{p, \theta}^{0, b}(\mathbb{T}^{m})$ and
$\mathbf{B}_{p, \theta}^{0, b}(\mathbb{T}^{m})$ are investigated in \cite{11}, \cite{12}, in particular, it is shown that these spaces do not coincide.

Recall the definition  spaces of functions with dominating mixed derivative $S_{p}^{\overline{r}}H,$ $S_{p, \theta}^{\overline{r}}B$ defined respectively by S. M.
Nikol'skii \cite{13} and T.I. Amanov (\cite{14} ch. I, sec. 17).

Let $\overline{r} = (r_{1},...,r_{m}),$ $\overline{k} = (k_{1},...,k_{m}),$ $k_{j}\in \mathbb{N}$, $k_{j}>r_{j} > 0$, $j =1,...,m,$
$1 \le p <\infty$, $1\leq \theta \le +\infty.$

Space $S_{p, \theta}^{\overline{r}}B$
--- consists of all functions $f\in {\mathring L}_{p} (\mathbb{T}^{m})$ for which
\begin{equation*}
\| f\|_{S_{p, \theta}^{\overline{r}}B
} = \|f\|_{p} + \Biggl[\int_{0}^{1}...\int_{0}^{1}\omega_{\bar k}^{\theta}(f, \overline{t})_{p}\prod_{j=1}^{m} \frac{dt_{j}}{t_{j}^{1+\theta r_{j}}}\Biggr]^{\frac{1}{\theta}} < +\infty.
\end{equation*}

P.\,I. ~Lizorkin and S.\,M.~Nikol'skii \cite{15} established equivalent norms of the space $S_{p, \theta}^{\overline{r}}B$ for $1 < p < +\infty,$ $1\le \theta \le \infty$ (in the case of $\theta=\infty$, see also \cite{16}, Theorem 4.4.6).
M. K. Potapov \cite{4}, \cite{5} defined and investigated the generalization of spaces
$S_{p}^{\overline{r}}H,$ $S_{p, \theta}^{\overline{r}}B$, with the replacement of the function $t^{r_{j}}$, for $r_{j} > 0$, $j=1,...,m$, by more general functions that satisfy certain conditions.

We define a space of functions with mixed logarithmic smoothness.
\begin{op}
Let $1\leq p< \infty$, $1\leq \tau < \infty$, $0< \theta \leq \infty$, $b_{j}> -1/\theta$, $j=1,...,m$. By $S_{p, \tau, \theta}^{0, \overline{b}}\mathbf{B}$ we denote the space of all functions $f\in {\mathring L}_{p, \tau} (\mathbb{T}^{m})$ for which
\begin{equation*}
\| f\|_{S_{p, \tau, \theta}^{0, \overline{b}}\mathbf{B}}
 = \|f\|_{p, \tau} + \Biggl[\int_{0}^{1}...\int_{0}^{1}\omega_{\overline{k}}^{\theta}(f, \overline{t})_{p, \tau}\prod_{j=1}^{m} \frac{(1-\log t_{j})^{\theta b_{j}}}{t_{j}}dt_{1}...dt_{m}\Biggr]^{\frac{1}{\theta}} < +\infty,
\end{equation*}
where $\overline{b}=(b_{1},...,b_{m})$, $\overline{k}=(k_{1},...,k_{m})$, $k_{j}\in \mathbb{N}$, $j=1,...,m$.  
\end{op}

We will also consider the following class $S_{p, \tau,\theta}^{0,\overline{b}}B$, consisting of all functions $f\in{\mathring L}_{p, \tau} (\mathbb{T}^{m})$ for which
 \begin{equation*}
\| f\|_{S_{p, \tau, \theta}^{0, \overline{b}}B} =
\left\{\sum\limits_{\bar s \in \Bbb{Z}_{+}^{m}} \prod_{j=1}^{m}(s_{j} + 1)^{b_{j}\theta} \|\delta_{\overline{s}}(f)\|_{p, \tau}^{\theta}\right\}^{\frac{1}{\theta}} < \infty,
 \end{equation*}
 for $1 < p < +\infty,$ $1\leq \tau< \infty$, $1\le \theta \le +\infty$, $\overline{b}=(b_{1},...,b_{m})$, $b_{j}> -1/\theta$, $j=1,...,m$.

In case when $\tau=p$, instead of $S_{p, p,\theta}^{0, \overline{b}}B$, $S_{p, p, \theta}^{0, \overline{b}}\mathbf{B}$ we write $S_{p, \theta}^{0, \overline{b}}B$, $S_{p, \theta}^{0, \overline{b}}\mathbf{B}$  respectively.

Note that the spaces $S_{p, \theta}^{0, \overline{b}}\mathbf{B}$ and $S_{p, \theta}^{0, \overline{b}}B$ are analogs of the spaces $B_{p, \theta}^{0, b}(\mathbb{T}^{m})$ and $\mathbf{B}_{p, \theta}^{0, b}$.

 The main goal of the paper is to find unimprovable conditions for embedding the spaces $S_{p, \theta}^{0, \overline{b}}\mathbf{B}$ and $S_{p, \theta}^{0, \overline{b }}B$ into each other for different ratios of the parameters of these spaces.
In addition to the introduction, the article consists of four sections. The first section contains auxiliary assertions necessary for proving the main results. 
 In the second section, we prove the embedding theorem for the space $S_{p, \tau_{1}, \theta}^{0, \overline{b}}B$ into ${\mathring L}_{p, \tau_{2}} (\mathbb{T}^{m})$, for $1< \tau_{2}<\tau_{1}< \infty$. The third section describes the spaces $S_{p, \theta}^{0, \overline{b}}\mathbf{B}$ and $S_{p, \theta}^{0, \overline{b}}B $ lacunar Fourier series in the trigonometric system.
 The main results are formulated and proved in the fourth section. Necessary and sufficient conditions for embedding spaces $S_{p, \tau, \theta}^{0, \overline{b}}\mathbf{B}$, $S_{p, \tau, \theta}^{0, \overline{b}}B$ are established here.

  \smallskip
\setcounter{equation}{0}
\setcounter{lemma}{0}
\setcounter{theorem}{0}

\section{Auxiliary statements}
\label{sec1}  

\smallskip
First, we introduce additional notation and give auxiliary statements. 
We will denote the set of indices $\{1,...,m\}$ by the symbol $e_{m},$ its arbitrary subset --- by $e$ and $|e|$ --- the number of elements of $e$ .

If given an element $\overline{r} =(r_{1},\ldots,r_{m})$, $m$ -- dimensional space with
non-negative coordinates, then $\overline{r}^{e} =(r_{1}^{e},\ldots,r_{m}^{e})$
vector with components $r_{j}^{e} = r_{j}$ for $j\in e$ and $r_{j}^{e}= 0$ for $j\notin e.$

\begin{theorem}\label{th1 1}  (see \cite[Theorem  1. 2]{17}).
Let $1 <p <\infty$, $1< \tau\leqslant 2$ or $2 <p <\infty$, $2 < \tau <\infty$. Then for any function $f\in L_{p,\tau}(\mathbb{T}^{m})$ the inequality holds
\begin{equation*}
\|f\|_{p, \tau} \leqslant C
\Bigl(\sum\limits_{\overline{s} \in \mathbb{Z}_{+}^{m}}\|\delta_{\overline{s}}\|_{p, \tau}^{\tau_{0}}\Bigr)^{1/\tau_{0}},
\end{equation*}
where $\tau_{0} = \min\{\tau, 2\}$.
\end{theorem}

\begin{theorem}\label{th1 2} (see \cite[Theorem 1. 3]{17}). Let $1 < p < \infty$, $1< \tau_{2} \leqslant 2$. 
If $\tau_{2} <\tau_{1}$ and the function $f\in L_{p,\tau_{1}}(\mathbb{T}^{m})$ satisfies the condition
\begin{equation*} 
\sum\limits_{\overline{s} \in \mathbb{Z}_{+}^{m}}\Bigl(\sum\limits_{j=1}^{m}(s_{j}+1)\Bigr)^{\tau_{2}(1/\tau_{2} - 1/\tau_{1})}\|\delta_{\overline{s}}\|_{p, \tau_{1}}^{\tau_{2}} < \infty,
\end{equation*}
then $f\in L_{p, \tau_{2}}(\mathbb{T}^{m})$ and the inequality holds
\begin{equation*} 
\|f\|_{p, \tau_{2}} \leqslant C\left(\sum\limits_{\overline{s}\in \mathbb{Z}_{+}^{m}}
\Bigl(\sum\limits_{j=1}^{m}(s_{j}+1)\Bigr)^{\tau_{2}(1/\tau_{2} - 1/\tau_{1})}\|\delta_{\overline{s}}\|_{p, \tau_{1}}^{\tau_{2}}\right)^{1/\tau_{2}}.
\end{equation*}
\end{theorem}

Further, we will often use the following well-known statement.
\begin{lemma}\label{lem1 1}
Let be given non-negative numbers $a_{\overline{s}}$, $b_{\overline{s}}$, for $\overline{s}\in\mathbb{Z}_{+}^{m}$ and $0< \beta, \theta \leq \infty$, $\alpha \in \mathbb{R}$. If   $0< \beta< \theta \leq \infty$, then
   \begin{equation*}
    \Bigl(\sum\limits_{\overline{s} \in \mathbb{Z}_{+}^{m}}a_{\overline{s}}^{\beta\alpha} b_{\overline{s}}^{\beta}\Bigr)^{\frac{1}{\beta}}<< \Bigl(\sum\limits_{\overline{s} \in \mathbb{Z}_{+}^{m}} b_{\overline{s}}^{\theta}\Bigr)^{\frac{1}{\theta}} \Bigl(\sum\limits_{\overline{s} \in \mathbb{Z}_{+}^{m}}a_{\overline{s}}^{\beta\alpha\eta^{`}}\Bigr)^{\frac{1}{\beta\eta^{`}}},
    \end{equation*}
 where  $\eta =\frac{\theta}{\beta}$, $\eta^{`}=\frac{\eta}{\eta - 1}$, if $\theta < \infty$ and $\eta^{`} =1$, if $\theta = \infty$.
 
If $0<  \theta \leq \beta < \infty$, then
  \begin{equation*}
    \Bigl(\sum\limits_{\overline{s} \in \mathbb{Z}_{+}^{m}}a_{\overline{s}}^{\beta\alpha} b_{\overline{s}}^{\beta}\Bigr)^{\frac{1}{\beta}} \leq \Bigl(\sum\limits_{\overline{s} \in \mathbb{Z}_{+}^{m}} b_{\overline{s}}^{\theta}\Bigr)^{\frac{1}{\theta}} \sup_{\overline{s} \in \mathbb{Z}_{+}^{m}}a_{\overline{s}}.
    \end{equation*}
 \end{lemma} 
This lemma is a consequence of the H\"{o}lder and Jensen inequality (see \cite{2}), respectively, in the cases of $0<\beta< \theta$ and $0< \theta \leq \beta < \infty$.

  \smallskip

\setcounter{equation}{0}
\setcounter{lemma}{0}
\setcounter{theorem}{0}

\section{Embedding theorems of spaces with mixed logarithmic smoothness in the Lorentz space}\label{sec2} 

\smallskip
  
\begin{theorem}\label{th2 1}  
Let $1 < p <+ \infty$ and $1< \tau_{2}\leq 2$ or $2 < p <+ \infty$ and $2< \tau_{2}< \infty$,  $0< \theta \leq \infty$, $\beta = \min\{2, \,\, \tau_{2}\}$ and numbers $b_{j}> - \frac{1}{\theta}$, for $j=1,\ldots, m$, $\tau_{2} < \tau_{1}< \infty$.

If  \begin{equation}\label{eq2 1}
    \sum\limits_{\overline{s} \in \mathbb{Z}_{+}^{m}}\Bigl(\sum\limits_{j=1}^{m}(s_{j}+1)\Bigr)^{\beta\eta^{'}(1/\tau_{2} - 1/\tau_{1})}\prod_{j=1}^{m}(s_{j} + 1)^{-\beta\eta^{'}b_{j}} < \infty,
\end{equation}
in the case $\beta < \theta < \infty$ and
 \begin{equation}\label{eq2 2}
 \sup_{\overline{s} \in \mathbb{Z}_{+}^{m}}\Bigl(\sum\limits_{j=1}^{m}(s_{j}+1)\Bigr)^{(1/\tau_{2} - 1/\tau_{1})}\prod_{j=1}^{m}(s_{j} + 1)^{-b_{j}}< \infty,  
\end{equation}  
in the case $\theta \leq\beta  < \infty$, then
$S_{p, \tau_{1}, \theta}^{0,  \overline{b}}B \subset L_{p, \tau_{2}}(\mathbb{T}^{m})$.
 \end{theorem}

 \proof
 Let $f \in S_{p, \tau_{1}, \theta}^{0,  \overline{b}}B$. It follows from Theorem 1. 1 and Theorem 1. 2 that
\begin{equation}\label{eq2 3}
\|f\|_{p, \tau_{2}} \ll \Bigl(\sum\limits_{\overline{s} \in \mathbb{Z}_{+}^{m}}\|\delta_{\overline{s}}(f)\|_{p, \tau_{2}}^{\beta}\Bigr)^{1/\beta} \ll \left(\sum\limits_{\overline{s}\in \mathbb{Z}_{+}^{m}}
\Bigl(\sum\limits_{j=1}^{m}(s_{j}+1)\Bigr)^{\beta(1/\tau_{2} - 1/\tau_{1})}\|\delta_{\overline{s}}\|_{p, \tau_{1}}^{\beta}\right)^{1/\beta}.
\end{equation}
Now applying Lemma 1. 1 when
\begin{equation*}
b_{\overline{s}}=\prod_{j=1}^{m}(s_{j} + 1)^{b_{j}}\|\delta_{\overline{s}}(f)\|_{p, \tau_{1}}, \, \, \, \, 
a_{\overline{s}} = \Bigl(\sum\limits_{j=1}^{m}(s_{j}+1)\Bigr)^{(1/\tau_{2} - 1/\tau_{1})}\prod_{j=1}^{m}(s_{j} + 1)^{-b_{j}}
\end{equation*}
for $\overline{s} \in \mathbb{Z}_{+}^{m}$ and according to  conditions \eqref{eq2 1} and \eqref{eq2 2} from \eqref{eq2 3}, we get $S_{p, \tau_{1}, \theta}^{0, \overline{b}}B\subset L_{p, \tau_{2}}(\mathbb{T}^{m})$.
   \hfill $\Box$

\smallskip

\setcounter{equation}{0}
\setcounter{lemma}{0}
\setcounter{theorem}{0}
 \setcounter{rem}{0}
 
\section{Description of the spaces $S_{p, \tau, \theta}^{0, \bar b}\mathbf{B}$ and $S_{p, \tau, \theta}^{0, \bar b}B$ by lacunary Fourier series}\label{sec3}  

\smallskip

\begin{op} 
Let $1 < p < +\infty,$ $1\leq \tau< \infty$. By $\Lambda_{p ,\tau}$ we denote the set of all functions $f\in {\mathring L}_{p, \tau} (\mathbb{T}^{m})$ that have a lacunary Fourier series
\begin{equation*} 
f(\overline{x}) \sim \sum_{\overline{\nu} \in \mathbb{Z}_{+}^{m}}\lambda_{\overline{\nu}} \prod_{j=1}^{m}\cos 2^{\nu_{j}}x_{j}, \,\, \, \, \lambda_{\overline{\nu}}\in \mathbb{R}.
\end{equation*}
\end{op}

In the case of $\tau=p$, the class $\Lambda_{p , p}=\Lambda_{p}$ is considered in \cite{6}.
\begin{lemma}\label{lem3 1}
Let $f \in \Lambda_{p ,\tau}$, $1 < p < +\infty,$ $1< \tau< \infty$.

1. If $1 < p < +\infty,$ $1\leq \tau< \infty$, then
\begin{equation*}
\|f\|_{p, \tau}\ll \Bigl(\sum_{\overline{\nu} \in \mathbb{Z}_{+}^{m}}|\lambda_{\overline{\nu}}|^{2}\Bigr)^{1/2}
\end{equation*}

2. If $1 < p \leq 2$ and $1< \tau \leq 2$ or  $2< p< +\infty,$ $1< \tau< \infty$,
then
\begin{equation*}
\Bigl(\sum_{\overline{\nu} \in \mathbb{Z}_{+}^{m}}|\lambda_{\overline{\nu}}|^{2}\Bigr)^{1/2}\ll \|f\|_{p, \tau}. 
\end{equation*}

3. If $1 < p \leq 2$ and $2< \tau < \infty$, then
\begin{equation*}
\Bigl(\sum_{\overline{\nu} \in \mathbb{Z}_{+}^{m}}|\lambda_{\overline{\nu}}|^{\tau}\Bigr)^{1/\tau}\ll \|f\|_{p, \tau}. 
\end{equation*}
\end{lemma}
 
 \proof
Since $\delta_{\overline{s}}(f, \overline{x}) = \lambda_{\overline{s}} \prod_{j=1}^{m}\cos 2^{s_{j}}x_{j}$, for the function $\Lambda_{p,\tau}$, then by the Littlewood--Paley theorem in Lorentz space \cite{16}, we have
\begin{equation*}
\|f\|_{p, \tau} \asymp \Bigl\|\Bigl(\sum_{\overline{s} \in \mathbb{Z}_{+}^{m}}\lambda_{\overline{s}}^{2}(\prod_{j=1}^{m}\cos 2^{s_{j}}x_{j})^{2} \Bigr)^{1/2}\Bigr\|_{p, \tau}\ll \Bigl(\sum_{\overline{s} \in \mathbb{Z}_{+}^{m}}\lambda_{\overline{s}}^{2}\Bigr)^{1/2},
\end{equation*}
for $1 < p < +\infty,$ $1< \tau< \infty$. The first statement is proved.

If $2< p< +\infty$, then $L_{p, \tau}(\mathbb{T}^{m}) \subset L_{2}(\mathbb{T}^{m})$  and $\|f\|_{2}\ll \|f\|_{p, \tau}$, for $1< \tau< \infty$. Therefore, taking into account Parseval's equality , we have
\begin{equation*}
\|f\|_{p, \tau}\gg \Bigl(\sum_{\overline{s} \in \mathbb{Z}_{+}^{m}}\lambda_{\overline{s}}^{2}\Bigr)^{1/2},
\end{equation*}
for $2< p< +\infty$, $1< \tau< \infty$.

If $1 <p\leq 2$ and $1<\tau\leq2$, then from Lemma 1. 5 \cite{16}, it follows that
\begin{equation*}
\|f\|_{p, \tau} \gg \Bigl(\sum\limits_{\overline{s} \in \mathbb{Z}_{+}^{m}}\|\delta_{\overline{s}}(f)\|_{p, \tau}^{2}\Bigr)^{1/2} \gg \Bigl(\sum_{\overline{s} \in \mathbb{Z}_{+}^{m}}\lambda_{\overline{s}}^{2}\Bigr)^{1/2}.
\end{equation*}

If $1 <p\leq 2$ and $2<\tau <\infty$, then by Lemma 1. 6 \cite{16}, we will get
\begin{equation*}
\|f\|_{p, \tau} \gg \Bigl(\sum\limits_{\overline{s} \in \mathbb{Z}_{+}^{m}}\|\delta_{\overline{s}}(f)\|_{p, \tau}^{\tau}\Bigr)^{1/\tau} \gg \Bigl(\sum_{\overline{s} \in \mathbb{Z}_{+}^{m}}|\lambda_{\overline{s}}|^{\tau}\Bigr)^{1/\tau}.
\end{equation*}
\hfill $\Box$

\begin{rem} In the case of $\tau = p$, Lemma 3.1 coincides with Lemma 3.9 \cite{6}.
\end{rem}

\begin{theorem}\label{th3 1}
Let $1 < p, \tau <+ \infty$, $0< \theta \leq \infty$ and the numbers $b_{j}>  \frac{1}{\theta}$, for $j=1,\ldots, m$.

1. If $f \in \Lambda_{p ,\tau}$, $1 < p < +\infty,$ $1\leq \tau< \infty$, then
\begin{equation*}
\|f\|_{S_{p, \tau, \theta}^{0, \overline{b}}\mathbf{B}} \ll \Biggl(\sum\limits_{\nu_{m} =1}^{\infty}...\sum\limits_{\nu_{1} =1}^{\infty}\prod_{j=1}^{m}(\nu_{j} + 1)^{\theta b_{j}}\Biggl(\sum\limits_{s_{m} =\nu_{m}}^{\infty}...\sum\limits_{s_{1} =\nu_{1}}^{\infty}|\lambda_{s_{1} - 1,...,s_{m} - 1}|^{2} \Biggr)^{\frac{\theta}{2}}
\Biggr)^{\frac{1}{\theta}}.
\end{equation*}

2. If $1 < p\leq 2$ and $1< \tau\leq 2$ or $2<p<+\infty$, $1< \tau< \infty$, then
for $f\in \Lambda_{p ,\tau}$ the inequality holds
\begin{equation*}
\Biggl(\sum\limits_{\nu_{m} =1}^{\infty}...\sum\limits_{\nu_{1} =1}^{\infty}\prod_{j=1}^{m}(\nu_{j} + 1)^{\theta b_{j}}\Biggl(\sum\limits_{s_{m} =\nu_{m}}^{\infty}...\sum\limits_{s_{1} =\nu_{1}}^{\infty}|\lambda_{s_{1} - 1,...,s_{m} - 1}|^{2} \Biggr)^{\frac{\theta}{2}}\Biggr)^{\frac{1}{\theta}} \ll\|f\|_{S_{p, \tau, \theta}^{0, \overline{b}}\mathbf{B}}.
\end{equation*}

3. If $1 < p \leq 2$ and $2< \to <\infty$, then the inequality holds
\begin{equation*}
\Biggl(\sum\limits_{\nu_{m} =1}^{\infty}...\sum\limits_{\nu_{1} =1}^{\infty}\prod_{j=1}^{m}(\nu_{j} + 1)^{\theta b_{j}}\Biggl(\sum\limits_{s_{m} =\nu_{m}}^{\infty}...\sum\limits_{s_{1} =\nu_{1}}^{\infty}|\lambda_{s_{1} - 1,...,s_{m} - 1}|^{\tau} \Biggr)^{\frac{\theta}{\tau}}\Biggr)^{\frac{1}{\theta}} \ll\|f\|_{S_{p, \tau, \theta}^{0, \overline{b}}\mathbf{B}}.
\end{equation*}
 for the function $f\in\Lambda_{p ,\tau}$
\end{theorem}

 \proof
 According to Lemma 3.1, the inequality holds for the function $f\in\Lambda_{p ,\tau}$
\begin{equation*}
\sigma_{\nu_{1},...,\nu_{m}}=\Biggl\|\Biggl(\sum\limits_{s_{m}=\nu_{m}}^{\infty}...\sum\limits_{s_{1}=\nu_{1}}^{\infty}|\delta_{\overline{s}}(f)|^{2}\Biggr)^{\frac{1}{2}}\Biggr\|_{p, \tau} \ll \Biggl(\sum\limits_{s_{m} =\nu_{m}}^{\infty}...\sum\limits_{s_{1} =\nu_{1}}^{\infty}|\lambda_{s_{1} - 1,...,s_{m} - 1}|^{2} \Biggr)^{\frac{1}{2}},
\end{equation*}
for $1 < p, \tau <+ \infty$. Therefore, by Theorem 2 \cite{19}, we will have
\begin{multline}\label{eq3 1}
\|f\|_{S_{p, \tau, \theta}^{0, \overline{b}}\mathbf{B}} \ll \Biggl\{\Biggl(\sum\limits_{s_{m} =1}^{\infty}...\sum\limits_{s_{1} =1}^{\infty}|\lambda_{s_{1} - 1,...,s_{m} - 1}|^{2} \Biggr)^{\frac{1}{2}}
\\
 + \Biggl(\sum\limits_{\nu_{m} =1}^{\infty}...\sum\limits_{\nu_{1} =1}^{\infty}\prod_{j=1}^{m}(\nu_{j} + 1)^{\theta b_{j}}\Biggl(\sum\limits_{s_{m} =\nu_{m}}^{\infty}...\sum\limits_{s_{1} =\nu_{1}}^{\infty}|\lambda_{s_{1} - 1,...,s_{m} - 1}|^{2} \Biggr)^{\frac{\theta}{2}}\Biggr)^{\frac{1}{\theta}} \Biggr\} .
\end{multline}
Since
\begin{multline*}
\sum\limits_{\nu_{m} =1}^{\infty}...\sum\limits_{\nu_{1} =1}^{\infty}\prod_{j=1}^{m}(\nu_{j} + 1)^{\theta b_{j}}\Biggl(\sum\limits_{s_{m} =\nu_{m}}^{\infty}...\sum\limits_{s_{1} =\nu_{1}}^{\infty}|\lambda_{s_{1} - 1,...,s_{m} - 1}|^{2} \Biggr)^{\frac{\theta}{2}} 
\\
\geq \prod_{j=1}^{m}2^{\theta b_{j}}\Biggl(\sum\limits_{s_{m} =1}^{\infty}...\sum\limits_{s_{1} =1}^{\infty}|\lambda_{s_{1} - 1,...,s_{m} - 1}|^{2} \Biggr)^{\frac{\theta}{2}},
\end{multline*}
then  the first statement follows from \eqref{eq3 1}.

The second and third statements follow from Lemma 3.1 and Theorem 2 \cite{19} in the corresponding values of the parameters $p, \, \, \tau$.
\hfill $\Box$

\begin{rem} In the case of $\tau = p$ from Theorem 3.1 we obtain an analogue of Theorem 5.3 \cite{12} for periodic functions.
\end{rem}

\begin{theorem}\label{th3 2} Let $1 < p, \tau <+ \infty$, $0< \theta \leq \infty$ and the numbers $b_{j}> - \frac{1}{\theta}$, for $j=1,\ldots, m$.
Then for the function $f\in\Lambda_{p ,\tau}$, the relation is valid
\begin{equation*}
\|f\|_{S_{p, \tau, \theta}^{0, \overline{b}}B} \asymp \Biggl(\sum\limits_{s_{m} =1}^{\infty}...\sum\limits_{s_{1} =1}^{\infty}\prod_{j=1}^{m}(s_{j} + 1)^{\theta b_{j}}|\lambda_{s_{1} - 1,...,s_{m} - 1}|^{\theta}\Biggr)^{\frac{1}{\theta}}.
\end{equation*}
\end{theorem}

 \proof
 For the function $f\in \Lambda_{p ,\tau}$, the relation is executed
\begin{equation*}
 \|\delta_{\overline{s}}(f)\|_{p, \tau} \asymp |\lambda_{s_{1} - 1,...,s_{m} - 1}|, \, \, \, s_{j}\in \mathbb{N}, \, \, j=1, \ldots. 
\end{equation*}
Therefore, the statement of the theorem follows from the definition of the space $S_{p, \tau, \theta}^{0, \overline{b}}B$.
 \hfill $\Box$
\begin{rem} In the case of $\tau = p$ from Theorem 3.2 we obtain an analogue of proposition 5.1 \cite{12} for periodic functions.
\end{rem}

\smallskip

\setcounter{equation}{0}
\setcounter{lemma}{0}
\setcounter{theorem}{0}
 \setcounter{rem}{0}
 
\section{On the relationship of function spaces with mixed logarithmic smoothness}\label{sec4} 
  
 \begin{theorem}\label{th4 1}
Let $1 < p <+ \infty$ and $1<\tau \leq 2$ or $2< p <+ \infty$ and $1< \tau  < \infty$, $0< \theta_{1}, \theta_{2} \leq \infty$ and the numbers $b_{j}^{(i)}> - \frac{1}{\theta_{i}}$, for $j=1,\ldots, m$, $\overline{b}^{(i)}=(b_{1}^{(i)}, \ldots , b_{m}^{(i)})$, $i=1, 2$. Then $S_{p, \tau, \theta_{1}}^{0, \overline{b}^{(1)}}\mathbf{B}\subset_{p, \tau, \theta_{2}}^{0, \overline{b}^{(2)}}\mathbf{B}$, if and only if

1.  $b_{j}^{(1)} + \frac{1}{\theta_{1}} > b_{j}^{(2)} + \frac{1}{\theta_{2}}$, $j=1,\ldots, m$;

2. $b_{j}^{(1)} + \frac{1}{\theta_{1}} = b_{j}^{(2)} + \frac{1}{\theta_{2}}$ and $0< \theta_{1} < \theta_{2} \leq \infty$.
\end{theorem}

\proof {\bf A sufficient part.}   
By Theorem 3 \cite{19} for the function $f\in S_{p, \tau, \theta_{i}}^{0, \overline{b}^{(i)}}\mathbf{B}$ the relation is valid
\begin{multline}\label{eq4 1} 
 \|f\|_{S_{p, \tau, \theta_{i}}^{0, \overline{b}}\mathbf{B}}
 \\
\asymp \Biggl\{ \|f\|_{p, \tau} + \Biggl(\sum\limits_{l_{m} =0}^{\infty}...\sum\limits_{l_{1} =0}^{\infty}\prod_{j=1}^{m}2^{l_{j}\theta (b_{j}^{(i)}+ \frac{1}{\theta_{i}})}\Biggl\|\sum\limits_{s_{m}=[2^{l_{m}-1}]+1}^{2^{l_{m}}}...\sum\limits_{s_{1}=[2^{l_{1}-1}]+1}^{2^{l_{1}}}\delta_{\overline{s}}(f)\Biggr\|_{p, \tau}^{\theta_{i}}\Biggr)^{\frac{1}{\theta_{i}}}\Biggr\}.    
\end{multline}
Let $f\in S_{p, \tau, \theta_{1}}^{0, \overline{b}^{(1)}}\mathbf{B}$. Then if
$0< \theta_{1} \leq \theta_{2} \leq\infty$, then according to Lemma 1. 1 we have
\begin{multline*}
\Biggl(\sum\limits_{l_{m} =0}^{\infty}...\sum\limits_{l_{1} =0}^{\infty}\prod_{j=1}^{m}2^{l_{j}\theta (b_{j}^{(2)}+ \frac{1}{\theta_{2}})}\Biggl\|\sum\limits_{s_{m}=[2^{l_{m}-1}]+1}^{2^{l_{m}}}...\sum\limits_{s_{1}=[2^{l_{1}-1}]+1}^{2^{l_{1}}}\delta_{\overline{s}}(f)\Biggr\|_{p, \tau}^{\theta_{2}}\Biggr)^{\frac{1}{\theta_{2}}}
\\
\ll \Biggl(\sum\limits_{l_{m} =0}^{\infty}...\sum\limits_{l_{1} =0}^{\infty}\prod_{j=1}^{m}2^{l_{j}\theta (b_{j}^{(1)}+ \frac{1}{\theta_{1}})}2^{\theta_{1}(b_{j}^{(1)}-b_{j}^{(2)} - \frac{1}{\theta_{1}} + \frac{1}{\theta_{2}})}\Biggl\|\sum\limits_{s_{m}=[2^{l_{m}-1}]+1}^{2^{l_{m}}}...\sum\limits_{s_{1}=[2^{l_{1}-1}]+1}^{2^{l_{1}}}\delta_{\overline{s}}(f)\Biggr\|_{p, \tau}^{\theta_{1}}\Biggr)^{\frac{1}{\theta_{1}}}.
\end{multline*}
Since in the first and second points of the number $b_{j}^{(1)} + \frac{1}{\theta_{1}}\geq  b_{j}^{(2)}+ \frac{1}{\theta_{2}}$, for $j=1,\ldots, m$, then from here we get $S_{p, \tau, \theta_{1}}^{0, \overline{b}^{(1)}}\mathbf{B}\subset S_{p, \tau, \theta_{2}}^{0, \overline{b}^{(2)}}\mathbf{B}$.

If $0< \theta_{2} < \theta_{1} \leq\infty$, then applying Lemma 1. 1 from the previous inequality we obtain the validity of the first statement.
 
{\bf The necessary part.}   Let the inclusion of  $S_{p, \tau, \theta_{1}}^{0, \overline{b}^{(1)}}\mathbf{B} \subset S_{p, \tau, \theta_{2}}^{0, \overline{b}^{(2)}}\mathbf{B}$ .

We assume that for some $j_{0}$ the inequality holds
$b_{j_{0}}^{(1)} + \frac{1}{\theta_{1}} < b_{j_{0}}^{(2)} + \frac{1}{\theta_{2}}$.
We will choose a number $\beta$ such that $\frac{1}{2} + b_{j_{0}}^{(1)} + \frac{1}{\theta_{1}} < \beta<b_{j_{0}}^{(2)} + \frac{1}{\theta_{2}} + \frac{1}{2}$. Let 's put
$\lambda_{s_{1}, ...,s_{m}} = (s_{j_{0}} + 1)^{-\beta}$, for $s_{j}=0$ for $j\neq j_{0}$ and $\lambda_{s_{1}, ...,s_{m}}= 0$, for $s_{j}=0$, $j=1, \ldots , m$.

 Now consider the function  
\begin{equation*}
f_{0}(\overline{x}) = \sum\limits_{s_{m} =0}^{\infty}...\sum\limits_{s_{1} =0}^{\infty}\lambda_{s_{1}, ...,s_{m}}\prod_{j=1}^{m}\cos 2^{s_{j}}x_{j} = \sum\limits_{s_{j_{0}} =0}^{\infty}\frac{\cos 2^{s_{j_{0}}}x_{j_{0}}}{(s_{j_{0}} + 1)^{\beta}}\prod_{j=1, j \neq j_{0}}^{m}\cos x_{j}.
\end{equation*}
Then using Theorem 3.1 we get
 \begin{multline}\label{eq4 2}
  \|f\|_{S_{p, \tau_{1}, \theta_{1}}^{0, \bar{b}^{(1)}}\mathbf{B}}
\asymp \Biggl(\sum\limits_{\nu_{j_{0}} =1}^{\infty}(\nu_{j_{0}} + 1)^{b_{j_{0}}^{(1)}\theta_{1}} \Biggl(\sum\limits_{s_{j_{0}} =\nu_{j_{0}}}^{\infty}(s_{j_{0}} + 1)^{-2\beta} \Biggr)^{\frac{\theta_{1}}{2}} \Biggr)^{\frac{1}{\theta_{1}}}
\\
\asymp  \Biggl(\sum\limits_{\nu_{j_{0}} =1}^{\infty} \frac{1}{(\nu_{j_{0}} + 1)^{(\beta -\frac{1}{2} -b_{j_{0}}^{(1)})\theta_{1}}}\Biggr)^{\frac{1}{\theta_{1}}}.
  \end{multline}

Since $(\beta -\frac{1}{2} -b_{j_{0}}^{(1)})\theta_{1} > 1$, then the series  
\begin{equation*}
  \sum\limits_{\nu_{j_{0}} =1}^{\infty} \frac{1}{(\nu_{j_{0}} + 1)^{(\beta -\frac{1}{2} -b_{j_{0}}^{(1)})\theta_{1}}}
  \end{equation*}
converges.
So it follows from \eqref{eq4 2} that the function $f_{0}\in S_{p, \tau, \theta_{1}}^{0, \overline{b}^{(1)}}\mathbf{B}$,   $1 < p <+ \infty$, $1< \tau  < \infty$, $0< \theta_{1} <\infty$.
 
Next, using the second statement of Theorem 3.1, we obtain
  \begin{multline}\label{eq4 3}
\|f\|_{S_{p, \tau, \theta_{2}}^{0, \bar{b}^{(2)}}\mathbf{B}}
\asymp \Biggl(\sum\limits_{\nu_{j_{0}} =1}^{\infty}(\nu_{j_{0}} + 1)^{b_{j_{0}}^{(2)}\theta_{2}} \Biggl(\sum\limits_{s_{j_{0}} =\nu_{j_{0}}}^{\infty}(s_{j_{0}} + 1)^{-2\beta} \Biggr)^{\frac{\theta_{2}}{2}} \Biggr)^{\frac{1}{\theta_{2}}}
\\
\asymp  \Biggl(\sum\limits_{\nu_{j_{0}} =1}^{\infty} \frac{1}{(\nu_{j_{0}} + 1)^{(\beta -\frac{1}{2} -b_{j_{0}}^{(2)})\theta_{2}}}\Biggr)^{\frac{1}{\theta_{2}}},
  \end{multline}
in the case $1< p\leq 2$, $1< \tau\leq 2$ or $2< p< \infty$, $2< \tau< \infty$.

Since by choosing the number $\beta$: $(\beta -\frac{1}{2} -b_{j_{0}}^{(2)})\theta_{2} < 1$, then the series  
\begin{equation*}
  \sum\limits_{\nu_{j_{0}} =1}^{\infty} \frac{1}{(\nu_{j_{0}} + 1)^{(\beta -\frac{1}{2} -b_{j_{0}}^{(2)})\theta_{2}}}
    \end{equation*}
diverges. 
Therefore, from the inequality \eqref{eq4 3} we get that $f_{0} \notin S_{p, \tau, \theta_{2}}^{0, \overline{b}^{(2)}}\mathbf{B}$, in the case of $1< p\leq 2$, $1< \tau_{2}\leq 2$ or $2< p< \infty$, $2< \tau< \infty$.

Thus, if there exists $j_{0}$ such that
$b_{j_{0}}^{(1)} + \frac{1}{\theta_{1}} < b_{j_{0}}^{(2)} + \frac{1}{\theta_{2}}$, then this contradicts the inclusion $S_{p, \tau, \theta_{1}}^{0, \overline{b}^{(1)}}\mathbf{B}\subset S_{p, \tau, \theta_{2}}^{0, \overline{b}^{(2)}}\mathbf{B}$, in the case of $1<p\leq 2$, $1< \tau\leq 2$ or $2< p< \infty$, $2< \tau < \infty$.  
 \hfill $\Box$
 
 \begin{rem} In the case of $\tau=p$, Theorem 4.1 is proved in \cite[Proposition 6. 3]{12}.
\end{rem} 
  
    \begin{theorem}\label{th4 2} 
  Let $1 < p <+ \infty$, $1< \tau_{2}< \tau_{1} < \infty$, $0< \theta_{1}, \theta_{2} \leq\infty$ and the numbers $b_{j}^{(i)}> - \frac{1}{\theta_{i}}$, for $j=1,\ldots, m$, $i = 1, 2$.
  
1. If $0< \theta_{2} < \theta_{1} \leq\infty$, then embedding $S_{p, \tau_{1}, \theta_{1}}^{0, \overline{b}^{(1)}}B \subset S_{p, \tau_{2}, \theta_{2}}^{0, \overline{b}^{(2)}}B$ is valid if and only if
\begin{equation}\label{eq4 4}
\Biggl(\sum\limits_{s_{m} =1}^{\infty}...\sum\limits_{s_{1} =1}^{\infty}\prod_{j=1}^{m}s_{j}^{(b_{j}^{(2)}-b_{j}^{(1)})\theta_{2}\eta^{'}}\Bigl(\sum\limits_{j =1}^{m}(s_{j} + 1)\Bigr)^{(\frac{1}{\tau_{2}}-\frac{1}{\tau_{1}})\theta_{2}\eta^{'}}\Biggr)^{\frac{1}{\theta_{2}\eta^{'}}} < \infty, 
\end{equation}
where 
 $\eta = \frac{\theta_{1}}{\theta_{2}}\, \, \,\,$, $\frac{1}{\eta} + \frac{1}{\eta^{`}} = 1$.

2. If $0< \theta_{1} \leq \theta_{2} \leq\infty$, then embedding $S_{p, \tau_{1}, \theta_{1}}^{0, \overline{b}^{(1)}}B\subset S_{p, \tau_{2}, \theta_{2}}^{0, \overline{b}^{(2)}}B$ is valid if 
\begin{equation}\label{eq4 5}
  \sup_{\overline{s}\in \mathbb{Z}_{+}^{m}}\prod_{j=1}^{m}(s_{j} +1)^{b_{j}^{(2)}-b_{j}^{(1)}}\Bigl(\sum\limits_{j =1}^{m}(s_{j} + 1)\Bigr)^{(\frac{1}{\tau_{2}}-\frac{1}{\tau_{1}})} < \infty.
\end{equation}
\end{theorem}
    
    \proof 
    Let
     $f \in S_{p, \tau_{1}, \theta_{1}}^{0, \overline{b}^{(1)}}B$. 
Since
$1< \tau_{2}< \tau_{1} < \infty$, then according to Theorem 1. 3 we have
\begin{multline}\label{eq4 6}
\Biggl(\sum\limits_{s_{m} =1}^{\infty}...\sum\limits_{s_{1} =1}^{\infty}\prod_{j=1}^{m}s_{j}^{b_{j}^{(2)}\theta_{2}}\|\delta_{\overline{s}}(f)\|_{p, \tau_{2}}^{\theta_{2}}\Biggr)^{\frac{1}{\theta_{2}}}
\\
\ll \Biggl(\sum\limits_{s_{m} =1}^{\infty}...\sum\limits_{s_{1} =1}^{\infty}\prod_{j=1}^{m}s_{j}^{b_{j}^{(2)}\theta_{2}}\Bigl(\sum\limits_{j =1}^{m}(s_{j} + 1)\Bigr)^{(\frac{1}{\tau_{2}}-\frac{1}{\tau_{1}})}\|\delta_{\overline{s}}(f)\|_{p, \tau_{1}}^{\theta_{2}}\Biggr)^{\frac{1}{\theta_{2}}}.
\end{multline}

The sufficient part of the first assertion is proved in \cite[Theorem 5]{19} 

If $0<\theta_{1}\leq\theta_{2}< \infty$,then according to Lemma 1. 1 and the condition\eqref{eq4 5}, from \eqref{eq4 6} it is not difficult to make sure that $S_{p, \tau_{1}, \theta_{1}}^{0, \overline{b}^{(1)}}B \subset S_{p, \tau_{2}, \theta_{2}}^{0, \overline{b}^{(2)}}B$. The second assertion is proven.

Now we prove the necessary part of the first assertion.
 Let $S_{p, \tau_{1}, \theta_{1}}^{0, \overline{b}^{(1)}}B\subset S_{p, \tau_{2}, \theta_{2}}^{0, \overline{b}^{(2)}}B$. Let's assume that the condition \eqref{eq4 5} is not met, i.e.
\begin{equation}\label{eq4 7}
\Biggl(\sum\limits_{s_{m} =1}^{\infty}...\sum\limits_{s_{1} =1}^{\infty}\prod_{j=1}^{m}s_{j}^{(b_{j}^{(2)}-b_{j}^{(1)})\theta_{2}\eta^{'}}\Bigl(\sum\limits_{j =1}^{m}(s_{j} + 1)\Bigr)^{(\frac{1}{\tau_{2}}-\frac{1}{\tau_{1}})\theta_{2}\eta^{'}}\Biggr)^{\frac{1}{\theta_{2}\eta^{'}}} = \infty. 
\end{equation}
Note that in the case of $\theta_{1} = \infty$, the number $\eta^{'}=1$.

We consider the set of all $\overline{s}\in\mathbb{Z}_{+}^{m}$ satisfying the inequalities $n\leq\sum\limits_{j=1}^{m}(s_{j}+1) <n+1$. Then \eqref{eq4 7} is equivalent to the following equality
\begin{equation}\label{eq4 8}
\sum\limits_{n =1}^{\infty}n^{(\frac{1}{\tau_{2}}-\frac{1}{\tau_{1}})\theta_{2}\eta^{'}}
\sum\limits_{n \leq \sum\limits_{j =1}^{m}(s_{j} + 1) < n+1}\prod_{j=1}^{m}s_{j}^{(b_{j}^{(2)}-b_{j}^{(1)})\theta_{2}\eta^{'}} = \infty. 
\end{equation}
For the brevity of the record , we introduce the notation:
$\delta = \theta_{2}\eta^{'}$,
  \begin{equation*}
\sigma_{n}:=\sum\limits_{n \leq \sum\limits_{j =1}^{m}(s_{j} + 1) < n+1}\prod_{j=1}^{m}s_{j}^{(b_{j}^{(2)}-b_{j}^{(1)})\theta_{2}\eta^{'}}, \, \, \, \, a_{n}:=n^{\frac{1}{\tau_{2}}-\frac{1}{\tau_{1}}}\sigma_{n}^{\frac{1}{\theta_{2}\eta^{'}}}.
\end{equation*}
Then according to \eqref{eq4 8}, according to the well-known Abel theorem (see \cite[p. 290]{20}), there exists a monotonically decreasing to zero numerical sequence $\{\varepsilon_{n}\}$ such that (see also \cite{8})
  \begin{equation*}
1) \, \, \,\,  \varepsilon_{n}a_{n}^{\delta} \leq 1; \, \, 2) \,\, \sum\limits_{n =1}^{\infty} \varepsilon_{n}^{\theta_{1}}a_{n}^{\delta} < \infty; \, \,\, \, 3) \, \,  
\sum\limits_{n =1}^{\infty} \varepsilon_{n}^{\theta_{2}}a_{n}^{\delta} = \infty.
\end{equation*}
We will introduce another notation 
 $B_{\overline{s}}=\prod_{j=1}^{m}s_{j}^{b_{j}^{(2)}-b_{j}^{(1)}}\Bigl(\sum\limits_{j =1}^{m}(s_{j} + 1)\Bigr)^{\frac{1}{\tau_{2}}-\frac{1}{\tau_{1}}}$, $\overline{s}\in \mathbb{Z}_{+}^{m}$

For $\overline{s}\in\mathbb{Z}_{+}^{m}$, we consider the function
\begin{equation*}
G_{\overline{s}}(\overline{x}) = \sum\limits_{j =1}^{m}\prod_{k \in \{1,...,m\}\setminus\{j\}}e^{i(2^{s_{k} - 1})2\pi x_{k}}\sum\limits_{\nu_{j} =2^{s_{j} - 1}}^{2^{s_{j}} - 1}\frac{\cos 2\pi\nu_{j}x_{j}}{\nu_{j}^{1-\frac{1}{p}}}, \, \, \, \, \overline{x}\in [0, 1)^{m}, 
\end{equation*}
 taking into account that $\eta^{'}=1$ for $\theta_{1}=\infty$.
In \cite{21} it is proved that
\begin{equation}\label{eq4 9}
\|G_{\overline{s}}\|_{p, \tau} \asymp \Bigl(\sum\limits_{j =1}^{m}(s_{j} + 1)\Bigr)^{\frac{1}{\tau}}, \, \, \, \, 1< p, \tau< \infty.  
\end{equation}
Now we will consider the function
\begin{equation*}
f_{3}(\overline{x}) = \sum_{n=1}^{\infty}\varepsilon_{n} \sum\limits_{n \leq \sum\limits_{j =1}^{m}(s_{j} + 1) < n+1}B_{\overline{s}}^{\frac{\delta}{\theta_{2}}}\prod_{j=1}^{m}(s_{j} + 1)^{-b_{j}^{(2)}}\Bigl(\sum\limits_{j =1}^{m}(s_{j} + 1)\Bigr)^{-\frac{1}{\tau_{2}}} G_{\overline{s}}(\overline{x}).
\end{equation*}
Then by definition of the space $S_{p, \tau_{1}, \theta_{1}}^{0, \overline{b}^{(1)}}B$ and taking into account the relation \eqref{eq4 9}, we will have
\begin{multline*} 
\|f_{3}\|_{S_{p, \tau_{1}, \theta_{1}}^{0, \overline{b}^{(1)}}B} = 
\Biggl(\sum\limits_{s_{m} =1}^{\infty}...\sum\limits_{s_{1} =1}^{\infty}\prod_{j=1}^{m}(s_{j}+1)^{b_{j}^{(1)}\theta_{1}}\|\delta_{\overline{s}}(f_{3})\|_{p, \tau_{1}}^{\theta_{1}}\Biggr)^{\frac{1}{\theta_{1}}}
\\
=\Biggl(\sum\limits_{n =1}^{\infty}\varepsilon_{n}^{\theta_{1}}\sum\limits_{n \leq \sum\limits_{j =1}^{m}(s_{j} + 1) < n+1} B_{\overline{s}}^{\frac{\delta}{\theta_{2}}\theta_{1}}\prod_{j=1}^{m}(s_{j} + 1)^{(b_{j}^{(1)}-b_{j}^{(2)})\theta_{1}}\Bigl(\sum\limits_{j =1}^{m}(s_{j} + 1)\Bigr)^{-\frac{\theta_{1}}{\tau_{2}}}\|G_{\overline{s}}\|_{p, \tau_{1}}^{\theta_{1}}\Biggr)^{\frac{1}{\theta_{1}}}
\\
\ll \Biggl(\sum\limits_{n =1}^{\infty}\varepsilon_{n}^{\theta_{1}}\sum\limits_{n \leq \sum\limits_{j =1}^{m}(s_{j} + 1) < n+1} B_{\overline{s}}^{\frac{\delta}{\theta_{2}}\theta_{1}}\prod_{j=1}^{m}(s_{j} + 1)^{(b_{j}^{(1)}-b_{j}^{(2)})\theta_{1}}\Bigl(\sum\limits_{j =1}^{m}(s_{j} + 1)\Bigr)^{(\frac{1}{\tau_{1}}-\frac{1}{\tau_{2}})\theta_{1}} \Biggr)^{\frac{1}{\theta_{1}}} 
\\
=C\Biggl(\sum\limits_{n =1}^{\infty}\varepsilon_{n}^{\theta_{1}}\sum\limits_{n \leq \sum\limits_{j =1}^{m}(s_{j} + 1) < n+1} B_{\overline{s}}^{\delta}\Biggr)^{\frac{1}{\theta_{1}}}  = C \Biggl(\sum\limits_{n =1}^{\infty}\varepsilon_{n}^{\theta_{1}} n^{\delta(\frac{1}{\tau_{2}}-\frac{1}{\tau_{1}})}\sigma_{n} \Biggr)^{\frac{1}{\theta_{1}}} = C \Biggl(\sum\limits_{n =1}^{\infty}\varepsilon_{n}^{\theta_{1}} a_{n}^{\delta} \Biggr)^{\frac{1}{\theta_{1}}}.
\end{multline*}
 
Therefore, according to the condition $2)$ sequence selection $\{\varepsilon_{n}\}$, it follows that the function $f_{3}\in S_{p, \tau_{1}, \theta_{1}}^{0, \overline{b}^{(1)}}B$, $1<p<\infty$, $1< \tau < \infty$, $0< \theta < \infty$.

Now we prove that the function $f_{3}\notin S_{p, \tau_{2}, \theta_{2}}^{0, \overline{b}^{(2)}}B$. Again, according to the relation \eqref{eq4 9}, we will have
\begin{multline*}
\|f_{3}\|_{S_{p, \tau_{2}, \theta_{2}}^{0, \overline{b}^{(2)}}B} \gg \Biggl(\sum\limits_{n =1}^{\infty}\varepsilon_{n}^{\theta_{2}}\sum\limits_{n \leq \sum\limits_{j =1}^{m}(s_{j} + 1) < n+1} \prod_{j=1}^{m}(s_{j} + 1)^{b_{j}^{(2)}\theta_{2}} B_{\overline{s}}^{\delta}\prod_{j=1}^{m}(s_{j} + 1)^{-b_{j}^{(2)}\theta_{2}}
\\
\times
\Bigl(\sum\limits_{j =1}^{m}(s_{j} + 1)\Bigr)^{-\frac{\theta_{2}}{\tau_{2}}}\|G_{\overline{s}}\|_{p, \tau_{2}}^{\theta_{2}}\Biggr)^{\frac{1}{\theta_{2}}}\gg \Biggl(\sum\limits_{n =1}^{\infty}\varepsilon_{n}^{\theta_{2}}\sum\limits_{n \leq \sum\limits_{j =1}^{m}(s_{j} + 1) < n+1} B_{\overline{s}}^{\delta}\Biggr)^{\frac{1}{\theta_{2}}} 
\\
\gg  \Biggl(\sum\limits_{n =1}^{\infty}\varepsilon_{n}^{\theta_{2}} n^{\delta(\frac{1}{\tau_{2}}-\frac{1}{\tau_{1}})} \sum\limits_{n \leq \sum\limits_{j =1}^{m}(s_{j} + 1) < n+1} \prod_{j=1}^{m}(s_{j} + 1)^{-(b_{j}^{(1)} - b_{j}^{(2)})\delta} \Biggr)^{\frac{1}{\theta_{2}}}= C \Biggl(\sum\limits_{n =1}^{\infty}\varepsilon_{n}^{\theta_{2}} a_{n}^{\delta} \Biggr)^{\frac{1}{\theta_{2}}}.
\end{multline*}
According to the condition $3)$ sequence selection $\{\varepsilon_{n}\}$, hence we get that the function $f_{3}\notin S_{p, \tau_{2}, \theta_{2}}^{0, \overline{b}^{(2)}}B$.
 \hfill $\Box$
 
\begin{theorem}\label{th4 3} 
Let $1 < p <+\infty$, $1< \tau_{2}< \tau_{1} < \infty$, $0< \theta_{1}, \theta_{2} \leq\infty$ and the numbers $b_{j}^{(i)}> - \frac{1}{\theta_{i}}$, for $j=1,\ldots, m$, $\overline{b}^{(i)}=(b_{1}^{(i)}, \ldots , b_{m}^{(i)})$, $i=1, 2$.

1.\, If $0< \theta_{2}< \theta_{1} \leq \infty$, then
$S_{p, \tau, \theta_{1}}^{0, \overline{b}^{(1)}}\mathbf{B}\subset S_{p, \tau, \theta_{2}}^{0, \overline{b}^{(2)}}\mathbf{B}$, if and only if 
\begin{equation}\label{eq4 10}
\sum\limits_{l_{m} =0}^{\infty}...\sum\limits_{l_{1} =0}^{\infty}\prod_{j=1}^{m} 2^{l_{j}(b_{j}^{(2)} - b_{j}^{(1)} - \frac{1}{\theta_{1}} + \frac{1}{\theta_{2}})\theta_{2}\eta^{`}}\Bigl(\sum_{j=1}^{m}2^{l_{j}}\Bigr)^{(\frac{1}{\tau_{2}} - \frac{1}{\tau_{1}})\theta_{2}\eta^{`}} < \infty,
\end{equation}
where $\eta = \frac{\theta_{1}}{\theta_{2}}$, $\frac{1}{\eta} + \frac{1}{\eta^{'}}$.

2. \,  If 
 $1< \tau_{2} < \tau_{1}  < \infty$,  $0< \theta_{1} < \theta_{2} \leq \infty$ and  
\begin{equation}\label{eq4 11}
 \sup_{\overline{s} \in \mathbb{Z}_{+}^{m}} \prod_{j=1}^{m}2^{s_{j}(b_{j}^{(2)} - b_{j}^{(1)} - \frac{1}{\theta_{1}} + \frac{1}{\theta_{2}})}\Bigl(\sum_{j=1}^{m}2^{s_{j}}\Bigr)^{\frac{1}{\tau_{2}} - \frac{1}{\tau_{1}}} < \infty,
\end{equation}  
then 
 $S_{p, \tau_{1}, \theta_{1}}^{0, \overline{b}^{(1)}}\mathbf{B} \subset S_{p, \tau_{2}, \theta_{2}}^{0, \overline{b}^{(2)}}\mathbf{B}$.
\end{theorem}

\proof We will first prove the first statement.  Sufficient part.  Let 
 $S_{p, \tau, \theta_{1}}^{0, \overline{b}^{(1)}}\mathbf{B} \subset S_{p, \tau, \theta_{2}}^{0, \overline{b}^{(2)}}\mathbf{B}$. 
Since $0< \theta_{2} < \theta_{1} \leq \infty$, then when $\eta=\frac{\theta_{1}}{\theta_{2}}$, $\frac{1}{\eta} + \frac{1}{\eta^{`}} = 1$ by Lemma 1. 1, we have
\begin{multline}\label{eq4 12}  
\Biggl(\sum\limits_{l_{m} =0}^{\infty}...\sum\limits_{l_{1} =0}^{\infty}\prod_{j=1}^{m}2^{l_{j}\theta_{2} (b_{j}^{(2)}+ \frac{1}{\theta_{2}})}\Bigl(\sum_{j=1}^{m}2^{l_{j}}\Bigr)^{(\frac{1}{\tau_{2}} - \frac{1}{\tau_{1}})\theta_{2}}
\\
\times \Biggl\|\sum\limits_{s_{m}=[2^{l_{m}-1}]+1}^{2^{l_{m}}}...\sum\limits_{s_{1}=[2^{l_{1}-1}]+1}^{2^{l_{1}}}\delta_{\overline{s}}(f)\Biggr\|_{p, \tau_{1}}^{\theta_{2}}\Biggr)^{\frac{1}{\theta_{2}}}
\\
\ll \Biggl(\sum\limits_{l_{m} =0}^{\infty}...\sum\limits_{l_{1} =0}^{\infty}\prod_{j=1}^{m}2^{l_{j}\theta_{1} (b_{j}^{(1)}+ \frac{1}{\theta_{1}})}\Biggl\|\sum\limits_{s_{m}=[2^{l_{m}-1}]+1}^{2^{l_{m}}}...\sum\limits_{s_{1}=[2^{l_{1}-1}]+1}^{2^{l_{1}}}\delta_{\overline{s}}(f)\Biggr\|_{p, \tau_{1}}^{\theta_{1}}\Biggr)^{\frac{1}{\theta_{1}}}
\\
\times \Biggl(\sum\limits_{l_{m} =0}^{\infty}...\sum\limits_{l_{1} =0}^{\infty}\prod_{j=1}^{m} 2^{l_{j}(b_{j}^{(2)} - b_{j}^{(1)} - \frac{1}{\theta_{1}} + \frac{1}{\theta_{2}})\theta_{2}\eta^{`}}\Bigl(\sum_{j=1}^{m}2^{l_{j}}\Bigr)^{(\frac{1}{\tau_{2}} - \frac{1}{\tau_{1}})\theta_{2}\eta^{`}} \Biggr)^{\frac{1}{\theta_{2}\eta^{`}}}
\end{multline}
If the condition \eqref{eq4 10} is satisfied,
then  from \eqref{eq4 1} and \eqref{eq4 12}, it follows that
\begin{multline}\label{eq4 13} 
\Biggl(\sum\limits_{l_{m} =0}^{\infty}...\sum\limits_{l_{1} =0}^{\infty}\prod_{j=1}^{m}2^{l_{j}\theta_{2} (b_{j}^{(2)}+ \frac{1}{\theta_{2}})}\Biggl\|\sum\limits_{s_{m}=[2^{l_{m}-1}]+1}^{2^{l_{m}}}...\sum\limits_{s_{1}=[2^{l_{1}-1}]+1}^{2^{l_{1}}}\delta_{\overline{s}}(f)\Biggr\|_{p, \tau_{2}}^{\theta_{2}}\Biggr)^{\frac{1}{\theta_{2}}}
\\
\ll \Biggl(\sum\limits_{l_{m} =0}^{\infty}...\sum\limits_{l_{1} =0}^{\infty}\prod_{j=1}^{m}2^{l_{j}\theta_{1} (b_{j}^{(1)}+ \frac{1}{\theta_{1}})}\Biggl\|\sum\limits_{s_{m}=[2^{l_{m}-1}]+1}^{2^{l_{m}}}...\sum\limits_{s_{1}=[2^{l_{1}-1}]+1}^{2^{l_{1}}}\delta_{\overline{s}}(f)\Biggr\|_{p, \tau_{1}}^{\theta_{1}}\Biggr)^{\frac{1}{\theta_{1}}}
\end{multline}
Now,  from the inequality \eqref{eq4 13}, it follows that $S_{p, \tau_{1}, \theta_{1}}^{0, \overline{b}^{(1)}}\mathbf{B} \subset S_{p, \tau_{2}, \theta_{2}}^{0, \overline{b}^{(2)}}\mathbf{B}$ in the case of $0< \theta_{2} < \theta_{1} \leq\infty$, provided \eqref{eq4 10}.

{\bf The necessary part.} 
Let $S_{p, \tau_{1}, \theta_{1}}^{0, \overline{b}^{(1)}}\mathbf{B} \subset S_{p, \tau_{2}, \theta_{2}}^{0, \overline{b}^{(2)}}\mathbf{B}$.
We assume that the condition \eqref{eq4 10} is not fulfilled, i.e. 
 \begin{equation}\label{eq4 14} 
\sum\limits_{l_{m} =0}^{\infty}...\sum\limits_{l_{1} =0}^{\infty}\prod_{j=1}^{m} 2^{l_{j}(b_{j}^{(2)} - b_{j}^{(1)} - \frac{1}{\theta_{1}} + \frac{1}{\theta_{2}})\theta_{2}\eta^{`}}\Bigl(\sum_{j=1}^{m}2^{l_{j}}\Bigr)^{(\frac{1}{\tau_{2}} - \frac{1}{\tau_{1}})\theta_{2}\eta^{`}} = \infty.
\end{equation}
We note that in the case of $\theta_{1} = \infty$, the number $\eta^{'}=1$.

We consider the set of all $\overline{\nu}\in\mathbb{Z}_{+}^{m}$ satisfying the inequalities $n\leq\sum\limits_{j=1}^{m}2^{\nu_{j}} <n+1$. Then \eqref{eq4 14} is equivalent to the following equality
\begin{equation}\label{eq4 15}  
\sum\limits_{n =1}^{\infty}n^{(\frac{1}{\tau_{2}}-\frac{1}{\tau_{1}})\theta_{2}\eta^{'}}
\sum\limits_{n \leq \sum\limits_{j =1}^{m}2^{\nu_{j}} < n+1}\prod_{j=1}^{m}2^{\nu_{j}(b_{j}^{(2)}-b_{j}^{(1)} -\frac{1}{\theta_{1}} + \frac{1}{\theta_{2}})\theta_{2}\eta^{'}} = \infty. 
\end{equation}
For the brevity of the record , we introduce the notation:  
$\delta = \theta_{2}\eta^{'}$,
  \begin{equation*}
\sigma_{n}:=\sum\limits_{n \leq \sum\limits_{j =1}^{m}2^{\nu_{j}} < n+1}\prod_{j=1}^{m}2^{\nu_{j}(b_{j}^{(2)}-b_{j}^{(1)} -\frac{1}{\theta_{1}} + \frac{1}{\theta_{2}})\theta_{2}\eta^{'}}, \, \, \, \, a_{n}:=n^{\frac{1}{\tau_{2}}-\frac{1}{\tau_{1}}}\sigma_{n}^{\frac{1}{\theta_{2}\eta^{'}}}.
\end{equation*}
Then, according to the assumption \eqref{eq4 15}, according to the well-known theorem of Abel \cite[p. 290]{20}, there exists a monotonically decreasing to zero numerical sequence $\{\varepsilon_{n}\}$, such that
  \begin{equation*}
1) \, \, \,\,  \varepsilon_{n}a_{n}^{\delta} \leq 1; \, \, 2) \,\, \sum\limits_{n =1}^{\infty} \varepsilon_{n}^{\theta_{1}}a_{n}^{\delta} < \infty; \, \,\, \, 3) \, \,  
\sum\limits_{n =1}^{\infty} \varepsilon_{n}^{\theta_{2}}a_{n}^{\delta} = \infty.
\end{equation*}
Let 's introduce another notation
 \begin{equation*}
 B_{\overline{s}}=\prod_{j=1}^{m}2^{s_{j}(b_{j}^{(2)}-b_{j}^{(1)}-\frac{1}{\theta_{1}} + \frac{1}{\theta_{2}})}\Bigl(\sum\limits_{j =1}^{m}2^{s_{j}}\Bigr)^{\frac{1}{\tau_{2}}-\frac{1}{\tau_{1}}}, \, \, \overline{s}\in \mathbb{Z}_{+}^{m}.
\end{equation*}
Now consider the function
\begin{equation*}
f_{4}(\overline{x}) = \sum_{n=1}^{\infty}\varepsilon_{n} \sum\limits_{n \leq \sum\limits_{j =1}^{m}2^{\nu_{j}} < n+1}B_{\overline{\nu}}^{\frac{\delta}{\theta_{2}}}\prod_{j=1}^{m}2^{-\nu_{j}(b_{j}^{(2)} + \frac{1}{\theta_{2}})}\Bigl(\sum\limits_{j =1}^{m}2^{\nu_{j}}\Bigr)^{-\frac{1}{\tau_{2}}} G_{\overline{\nu}}(\overline{x}),
\end{equation*}
where  
\begin{equation*}
G_{\overline{\nu}}(\overline{x}) = \sum\limits_{j =1}^{m}\prod_{k \in \{1,...,m\}\setminus\{j\}}e^{i(2^{2^{\nu_{k}}} - 1)2\pi x_{k}}\sum\limits_{s_{j} =2^{\nu_{j}} + 1}^{2^{\nu_{j} + 1}}\sum\limits_{l_{j} =2^{s_{j} - 1}}^{2^{s_{j}} - 1}\frac{\cos 2\pi l_{j}x_{j}}{l_{j}^{1-\frac{1}{p}}}, \, \, \, \, \overline{x}\in [0, 1)^{m}, 
\end{equation*}
for  $\overline{s}\in \mathbb{Z}_{+}^{m}$.
It is known that \cite{21}
\begin{equation}\label{eq4 16} 
\|G_{\overline{\nu}}\|_{p, \tau} \asymp \Bigl(\sum\limits_{j =1}^{m}2^{\nu_{j}}\Bigr)^{\frac{1}{\tau}}, \, \, \, \, 1< p, \tau< \infty.  
\end{equation}
We prove that the function $f_{4} \in S_{p, \tau_{1}, \theta_{1}}^{0, \overline{b}^{(1)}}\mathbf{B}$.
According to the relation \eqref{eq4 16}, we will have
\begin{multline*} 
\sum\limits_{\overline{\nu} \in \mathbb{Z}_{+}^{m}}\prod_{j=1}^{m}2^{\nu_{j}(b_{j}^{(1)} + \frac{1}{\theta_{1}})\theta_{1}}\Biggl\|\sum\limits_{s_{1} =[2^{\nu_{1} - 1}] + 1}^{2^{\nu_{1}}} \ldots \sum\limits_{s_{m} =[2^{\nu_{m} - 1}] + 1}^{2^{\nu_{m}}}\delta_{\overline{s}}(f_{4})\Biggr\|_{p, \tau_{1}}^{\theta_{1}}
\\
=\sum\limits_{n =1}^{\infty}\varepsilon_{n}^{\theta_{1}}\sum\limits_{n \leq \sum\limits_{j =1}^{m}2^{\nu_{j}} < n+1}\prod_{j=1}^{m}2^{\nu_{j}(b_{j}^{(1)} + \frac{1}{\theta_{1}})\theta_{1}} B_{\overline{\nu}}^{\frac{\delta}{\theta_{2}}\theta_{1}}\prod_{j=1}^{m}2^{-\nu_{j}(b_{j}^{(2)} + \frac{1}{\theta_{2}})\theta_{1}}\Bigl(\sum\limits_{j =1}^{m}2^{\nu_{j}}\Bigr)^{-\frac{\theta_{1}}{\tau_{2}}}\|G_{\overline{\nu}}\|_{p, \tau_{1}}^{\theta_{1}}
\end{multline*}
\begin{multline*}
\ll \sum\limits_{n =1}^{\infty}\varepsilon_{n}^{\theta_{1}}\sum\limits_{n \leq \sum\limits_{j =1}^{m}2^{\nu_{j}} < n+1} B_{\overline{\nu}}^{\frac{\delta}{\theta_{2}}\theta_{1}}\prod_{j=1}^{m}2^{\nu_{j}(b_{j}^{(1)}-b_{j}^{(2)} - \frac{1}{\theta_{2}} + \frac{1}{\theta_{1}})\theta_{1}}\Bigl(\sum\limits_{j =1}^{m}2^{\nu_{j}}  \Bigr)^{(\frac{1}{\tau_{1}}-\frac{1}{\tau_{2}})\theta_{1}}  
\\
=C\sum\limits_{n =1}^{\infty}\varepsilon_{n}^{\theta_{1}}\sum\limits_{n \leq \sum\limits_{j =1}^{m}2^{\nu_{j}} < n+1}B_{\overline{\nu}}^{\frac{\delta}{\theta_{2}}\theta_{1}}B_{\overline{\nu}}^{-\theta_{1}} =
C\sum\limits_{n =1}^{\infty}\varepsilon_{n}^{\theta_{1}}\sum\limits_{n \leq \sum\limits_{j =1}^{m}2^{\nu_{j}} < n+1} B_{\overline{\nu}}^{\delta} 
\\
 = C \sum\limits_{n =1}^{\infty}\varepsilon_{n}^{\theta_{1}} n^{\delta(\frac{1}{\tau_{2}}-\frac{1}{\tau_{1}})}\sigma_{n}  = C \sum\limits_{n =1}^{\infty}\varepsilon_{n}^{\theta_{1}} a_{n}^{\delta}.
\end{multline*}
Therefore, according to the condition $2)$ sequence selection $\{\varepsilon_{n}\}$, it follows that the function $f_{4}\in S_{p, \tau_{1}, \theta_{1}}^{0, \overline{b}^{(1)}}\mathbf{B}$, $1< p< \infty$, $1< \tau < \infty$, $0< \theta_{1}< \infty$.

Now we prove that the function $f_{4} \notin S_{p, \tau_{2}, \theta_{2}}^{0, \overline{b}^{(2)}}\mathbf{B}$. Again taking into account the ratio \eqref{eq4 16} we will have
\begin{multline*}
\|f_{4}\|_{S_{p, \tau_{2}, \theta_{2}}^{0, \overline{b}^{(2)}}\mathbf{B}} \gg \Biggl(\sum\limits_{n =1}^{\infty}\varepsilon_{n}^{\theta_{2}}\sum\limits_{n \leq \sum\limits_{j =1}^{m}2^{\nu_{j}} < n+1} \prod_{j=1}^{m}2^{\nu_{j}(b_{j}^{(2)}+ \frac{1}{\theta_{2}})\theta_{2}} B_{\overline{\nu}}^{\delta}\prod_{j=1}^{m}2^{-\nu_{j}(b_{j}^{(2)} + \frac{1}{\theta_{2}} )\theta_{2}}
\\
\times
\Bigl(\sum\limits_{j =1}^{m}2^{\nu_{j}}\Bigr)^{-\frac{\theta_{2}}{\tau_{2}}}\|G_{\overline{\nu}}\|_{p, \tau_{2}}^{\theta_{2}}\Biggr)^{\frac{1}{\theta_{2}}}\gg \Biggl(\sum\limits_{n =1}^{\infty}\varepsilon_{n}^{\theta_{2}}\sum\limits_{n \leq \sum\limits_{j =1}^{m}2^{\nu_{j}} < n+1} B_{\overline{\nu}}^{\delta}\Biggr)^{\frac{1}{\theta_{2}}} 
\\
\gg  \Biggl(\sum\limits_{n =1}^{\infty}\varepsilon_{n}^{\theta_{2}} n^{\delta(\frac{1}{\tau_{2}}-\frac{1}{\tau_{1}})} \sum\limits_{n \leq \sum\limits_{j =1}^{m}2^{\nu_{j}} < n+1} \prod_{j=1}^{m}2^{-\nu_{j}(b_{j}^{(1)} - b_{j}^{(2)} -\frac{1}{\theta_{2}} + \frac{1}{\theta_{1}})\delta} \Biggr)^{\frac{1}{\theta_{2}}}= C \Biggl(\sum\limits_{n =1}^{\infty}\varepsilon_{n}^{\theta_{2}} a_{n}^{\delta} \Biggr)^{\frac{1}{\theta_{2}}}.
\end{multline*}
By condition $3)$ sequence selection $\{\varepsilon_{n}\}$, hence we get the function $f_{4} \notin S_{p, \tau_{2}, \theta_{2}}^{0, \overline{b}^{(2)}}\mathbf{B}$.
 
 Let's prove the second statement. Since $1< \tau_{2} < \tau_{1} < \infty$, then according to Theorem 1. 3 and Lemma 1. 1 for $0< \theta_{1} \leq\theta_{2} < \infty$, we have
\begin{multline*}
\Biggl(\sum\limits_{l_{m} =0}^{\infty}...\sum\limits_{l_{1} =0}^{\infty}\prod_{j=1}^{m}2^{l_{j}\theta_{2} (b_{j}^{(2)}+ \frac{1}{\theta_{2}})}\Biggl\|\sum\limits_{s_{m}=[2^{l_{m}-1}]+1}^{2^{l_{m}}}...\sum\limits_{s_{1}=[2^{l_{1}-1}]+1}^{2^{l_{1}}}\delta_{\overline{s}}(f)\Biggr\|_{p, \tau_{2}}^{\theta_{2}}\Biggr)^{\frac{1}{\theta_{2}}}
\\
\ll \Biggl(\sum\limits_{l_{m} =0}^{\infty}...\sum\limits_{l_{1} =0}^{\infty}\prod_{j=1}^{m}2^{l_{j}\theta_{2} (b_{j}^{(2)}+ \frac{1}{\theta_{2}})}\Bigl(\sum_{j=1}^{m}2^{l_{j}}\Bigr)^{(\frac{1}{\tau_{2}} - \frac{1}{\tau_{1}})\theta_{2}}\Biggl\|\sum\limits_{s_{m}=[2^{l_{m}-1}]+1}^{2^{l_{m}}}...\sum\limits_{s_{1}=[2^{l_{1}-1}]+1}^{2^{l_{1}}}\delta_{\overline{s}}(f)\Biggr\|_{p, \tau_{1}}^{\theta_{2}}\Biggr)^{\frac{1}{\theta_{2}}}
\end{multline*}
\begin{multline}\label{eq4 17}
\ll  \Biggl(\sum\limits_{l_{m} =0}^{\infty}...\sum\limits_{l_{1} =0}^{\infty}\prod_{j=1}^{m}2^{l_{j}\theta_{1} (b_{j}^{(1)}+ \frac{1}{\theta_{1}})}\Biggl\|\sum\limits_{s_{m}=[2^{l_{m}-1}]+1}^{2^{l_{m}}}...\sum\limits_{s_{1}=[2^{l_{1}-1}]+1}^{2^{l_{1}}}\delta_{\overline{s}}(f)\Biggr\|_{p, \tau_{1}}^{\theta_{1}}\Biggr)^{\frac{1}{\theta_{1}}} 
\\
\times
\sup_{\overline{s} \in \mathbb{Z}_{+}^{m}} \prod_{j=1}^{m}2^{s_{j}(b_{j}^{(2)} - b_{j}^{(1)} - \frac{1}{\theta_{1}} + \frac{1}{\theta_{2}})}\Bigl(\sum_{j=1}^{m}2^{s_{j}}\Bigr)^{\frac{1}{\tau_{2}} - \frac{1}{\tau_{1}}}.    
\end{multline}
From the relation \eqref{eq4 1} and the inequality \eqref{eq4 17} according to the condition \eqref{eq4 11} we get $S_{p, \tau_{1}, \theta_{1}}^{0, \overline{b}^{(1)}}\mathbf{B} \subset S_{p, \tau_{2}, \theta_{2}}^{0, \overline{b}^{(2)}}\mathbf{B}$.
 \hfill $\Box$

  \begin{theorem}\label{th4 4} 
Let $1 <p, \to <+ \infty$, $0< \theta \leq \infty$ and the numbers $b_{j}> - \frac{1}{\theta}$, for $j=1,\ldots, m$, $\min\{2, \tau, \theta\}=2$. For any number $\varepsilon > 0$, there is a function $f\in S_{p, \tau, \theta}^{0, \overline{v} - \varepsilon\overline{e}_{j_{0}}}B$, such that $f\notin S_{p, \tau, \theta}^{0, \overline{b}}\mathbf{B}$, where $\overline{e}_{j_{0}}$ is a point with coordinates $x_{j}=0$, for $j\neq j_{0}$ and $x_{j_{0}}=1$, $\overline{v}=(v_{1}, \ldots ,v_{m})$,
$v_{j}=b_{j}+\frac{1}{\min\{\beta, \theta\}}$,
  \begin{equation*}
\beta =\left\{\begin{array}{rl} \tau, & \mbox{if} \, \,  1<\tau \leq 2 \, \, \mbox{and} \, \, 1<p< \infty \\
 2, & \mbox{if} \, \, 2< \tau < \infty \, \,  \mbox{and} \, \, 2<p< \infty ;
\end{array}
 \right.    
   \end{equation*} 
\end{theorem}
    
 \proof 
 We will choose the number $\mu\in (\frac{1}{\theta}, \frac{1}{2})$.
We will consider the function   
\begin{equation*}
f_{4}(\overline{x}) = \sum\limits_{s_{j_{0}} =1}^{\infty}(s_{j_{0}} +1)^{-(v_{j_{0}} - \varepsilon + \frac{1}{\theta})}(1 +\log s_{j_{0}})^{-\mu}\cos 2^{s_{j_{0}}}x_{j_{0}}\prod_{j\neq j_{0}}\cos x_{j} .
\end{equation*}
Since $\mu\theta > 1$, then for this function, the series
\begin{equation*}
    \sum\limits_{s_{m} =1}^{\infty}...\sum\limits_{s_{1} =1}^{\infty}\prod_{j=1}^{m}s_{j}^{(v_{j} - \varepsilon\overline{e}_{j_{0}})\theta}|\lambda_{s_{1},...,s_{m}}|^{\theta} =\sum\limits_{s_{j_{0}} =1}^{\infty}\frac{1}{(s_{j_{0}} +1)(1 +\log s_{j_{0}})^{\mu\theta}}
\end{equation*}
converges.
Therefore, by Theorem 3.2, the function $f_{4}\in S_{p, \tau, \theta}^{0, \overline{v} - \varepsilon\overline{e}_{j_{0}}}B$ .
 
 We show that if $\min\{\beta, \theta\} = 2$, then the function $f_{4} \notin S_{p, \tau, \theta}^{0, \overline{b}}\mathbf{B}$.
 
 Since for any number $\varepsilon > 0$, the sequence $\{(s_{j_{0}} +1)^{\varepsilon}(1 +\log s_{j_{0}})^{-\mu}\}$ increases and $(v_{j_{0}} + \frac{1}{\theta})2 > $ 1, then  
\begin{multline*} 
\sum\limits_{s_{j_{0}} =\nu_{j_{0}}}^{\infty}((s_{j_{0}} +1)^{-(v_{j_{0}} - \varepsilon + \frac{1}{\theta})}(1 +\log s_{j_{0}})^{-\mu})^{2}
\\
 \gg (\nu_{j_{0}} +1)^{2\varepsilon}(1 +\log \nu_{j_{0}})^{-2\mu}  \sum\limits_{s_{j_{0}} =\nu_{j_{0}}}^{\infty} (s_{j_{0}} +1)^{-(v_{j_{0}} + \frac{1}{\theta})2} 
\\ 
 \gg  (\nu_{j_{0}} +1)^{2\varepsilon}(1 +\log \nu_{j_{0}})^{-2\mu}(\nu_{j_{0}} +1)^{1-(v_{j_{0}} + \frac{1}{\theta})2} = (\nu_{j_{0}} +1)^{2\varepsilon}(1 +\log \nu_{j_{0}})^{-2\mu}(\nu_{j_{0}} +1)^{-(b_{j_{0}} + \frac{1}{\theta})2}.
\end{multline*}
Therefore, taking into account that the sequence $\{(\nu_{j_{0}} +1)^{\varepsilon}(1 +\log \nu_{j_{0}})^{\mu(2-\theta)}\}$ increases we will have
    \begin{multline}\label{eq4 18}
\sum\limits_{\nu_{j_{0}}=1}^{\infty}(\nu_{j_{0}} +1)^{b_{j_{0}\theta}}\Biggl(\sum\limits_{s_{j_{0}} =\nu_{j_{0}}}^{\infty}((s_{j_{0}} +1)^{-(v_{j_{0}} - \varepsilon + \frac{1}{\theta})}(1 +\log s_{j_{0}})^{-\mu})^{2} \Biggr)^{\theta/2} 
\\
=\sum\limits_{\nu_{j_{0}}=1}^{\infty} \frac{(\nu_{j_{0}} +1)^{\varepsilon\theta}}{(\nu_{j_{0}} +1)(1 +\log \nu_{j_{0}})^{\mu\theta}} =\sum\limits_{\nu_{j_{0}}=1}^{\infty} \frac{(\nu_{j_{0}} +1)^{\varepsilon\theta}(1 +\log \nu_{j_{0}})^{\mu(2-\theta)}}{(\nu_{j_{0}} +1)(1 +\log \nu_{j_{0}})^{2\mu}} 
\\
\gg \sum\limits_{\nu_{j_{0}}=1}^{\infty} \frac{1}{(\nu_{j_{0}} +1)(1 +\log \nu_{j_{0}})^{2\mu}}.
\end{multline}
Since by choice $2\mu < 1$, then from \eqref{eq4 15}, it follows that for the function $f_{4}$ the series
 \begin{equation*}
    \sum\limits_{\nu_{m} =1}^{\infty}...\sum\limits_{\nu_{1} =1}^{\infty}\prod_{j=1}^{m}(\nu_{j}+1)^{b_{j}\theta} \Biggl(\sum\limits_{s_{1} =\nu_{1}}^{\infty}...\sum\limits_{s_{m} =\nu_{m}}^{\infty}  |\lambda_{s_{1},...,s_{m}}|^{2}\Biggr)^{\theta/2}
\end{equation*}
  diverges.
Therefore, by Theorem 3.1, the function $f_{4}\notin S_{p, \tau, \theta}^{0, \overline{b}}\mathbf{B}$, in the case of $2<p<\infty$ and $2<\tau < \infty$.  
 \hfill $\Box$    
     
\begin{theorem}\label{th4 5} 
Let $1 < p, \tau <+ \infty$, $0< \theta \leq \infty$ and the numbers $b_{j}> - \frac{1}{\theta}$,   $u_{j}=b_{j}+\frac{1}{\max\{\gamma, \theta\}}$, for $j=1,\ldots, m$,  $\overline{u}=(u_{1}, \ldots ,u_{m})$, where
\begin{equation*}
\gamma =\left\{\begin{array}{rl} 2, & \mbox{if} \, \,  1<\tau \leq 2 \, \, \mbox{and} \, \, 1<p< \infty \\
 \tau, & \mbox{if} \, \, 2< \tau < \infty \, \,  \mbox{and} \, \, 2<p< \infty .
\end{array}
 \right.    
   \end{equation*}
   If $\max\{\gamma, \theta\}=2$, then for any number $\varepsilon >0$ , the space $S_{p, \tau, \theta}^{0, \overline{b}}\mathbf{B}$ will not be nested in $S_{p, \tau, \theta}^{0, \overline{u} + \varepsilon\overline{e}_{j_{0}}}B$, where $\overline{e}_{j_{0}}$ is a point of space $\mathbb{R}^{m}$, with coordinates $x_{j}=0$, for $j\neq j_{0}$ and $x_{j_{0}}=1$.
   \end{theorem}
\proof 
For a given number $\varepsilon > 0$, choose a number $\mu$ such that $b_{j_{0}} + \frac{1}{\theta} + \frac{1}{\max\{\gamma, \theta\}} < \mu \leq b_{j_{0}} + \frac{1}{\theta} + \frac{1}{\max\{\gamma, \theta\}} + \varepsilon$.  

We will consider the function
\begin{equation*}
f_{5}(\overline{x}) = \sum\limits_{s_{j_{0}} =1}^{\infty}(s_{j_{0}} +1)^{-\mu}\cos 2^{s_{j_{0}}}x_{j_{0}}\prod_{j\neq j_{0}}\cos x_{j} .
\end{equation*}
Since $2\mu > 2(b_{j_{0}} + \frac{1}{\theta} + \frac{1}{2})> 1$, then for the function $f_{5}$ we will have 
 \begin{multline*}
      \sum\limits_{\nu_{m} =1}^{\infty}...\sum\limits_{\nu_{1} =1}^{\infty}\prod_{j=1}^{m}(\nu_{j}+1)^{b_{j}\theta} \Biggl(\sum\limits_{s_{1} =\nu_{1}}^{\infty}...\sum\limits_{s_{m} =\nu_{m}}^{\infty}  |\lambda_{s_{1},...,s_{m}}|^{2}\Biggr)^{\theta/2}
    \\
    = \sum\limits_{\nu_{j_{0}}=1}^{\infty}(\nu_{j_{0}} +1)^{b_{j_{0}}\theta}\Biggl(\sum\limits_{s_{j_{0}} =\nu_{j_{0}}}^{\infty}(s_{j_{0}} +1)^{-\mu 2} \Biggr)^{\theta/2} \ll \sum\limits_{\nu_{j_{0}}=1}^{\infty} \frac{1}{(\nu_{j_{0}} +1)^{\theta(\mu - \frac{1}{2} - b_{j_{0}})}} < \infty.
\end{multline*}
Therefore, by Theorem 3.1, the function $f_{5}\in S_{p, \tau, \theta}^{0, \overline{b}}\mathbf{B}$, $1 <p, \tau<+ \infty$, $0<\theta < \infty$.

Now we prove that $f_{5} \notin S_{p, \tau, \theta}^{0, \overline{u} + \varepsilon\overline{e}_{j_{0}}}B$. Since $\mu\leq b_{j_{0}} + \frac{1}{\theta} + \frac{1}{\max\{\gamma, \theta\}} + \varepsilon$, then $\theta(\mu - u_{j_{0}} - \varepsilon) \leq 1$. Therefore, for the function $f_{5}$ we have
\begin{equation*}
    \sum\limits_{s_{m} =1}^{\infty}...\sum\limits_{s_{1} =1}^{\infty}\prod_{j=1}^{m}s_{j}^{u_{j}\theta}|\lambda_{s_{1},...,s_{m}}|^{\theta} =\sum\limits_{s_{j_{0}} =1}^{\infty}\frac{1}{(s_{j_{0}} +1)^{\theta(\beta - u_{j_{0}} -\varepsilon)}} = \infty.
\end{equation*}
Therefore, according to Theorem 3.2, the function $f_{5} \notin S_{p, \tau, \theta}^{0, \overline{u} + \varepsilon\overline{e}_{j_{0}}}B$.
 \hfill $\Box$
 
 \begin{rem} In the case of $\tau = p$, Theorem 4.3 and Theorem 4.4 are analogs of statements 9. 1 and 9. 2 in [12].
\end{rem}

\end{document}